\documentclass[journal,twoside,web]{ieeecolor}
\usepackage{lcsys}
\usepackage{cite}
\usepackage{amsmath,amssymb,amsfonts}
\usepackage{algorithmic}
\usepackage{graphicx}
\usepackage{textcomp}

\newtheorem{theorem}{Theorem}
\newtheorem{defi}{Definition}
\newtheorem{remark}{Remark}
\newtheorem{assump}{Assumption}

\newcommand*{\QEDB}{\hfill\ensuremath{\square}}%

\def\BibTeX{{\rm B\kern-.05em{\sc i\kern-.025em b}\kern-.08em
    T\kern-.1667em\lower.7ex\hbox{E}\kern-.125emX}}
\begin{document}
\title{Recurrent Neural Network ODE Output for Classification Problems Follows the Replicator Dynamics}
\author{Julian Barreiro-Gomez and Jorge I. Poveda
\thanks{J. Barreiro-Gomez was supported in part by the Center on Stability, Instability, and Turbulence, and  Tamkeen under the NYU Abu Dhabi Research Institute grant CG002. J. I. Poveda was supported in part by the NSF grant ECCS 2305756.}
\thanks{J. Barreiro-Gomez is with NYU AD Research Institute, New York University Abu Dhabi, PO Box 129188, UAE (e-mail: jbarreiro@nyu.edu). }
\thanks{J. I. Poveda  is  with  the ECE Department, University of California, San Diego, CA USA. (e-mail:  poveda@ucsd.edu).}}

\maketitle
\thispagestyle{empty}

\begin{abstract}
This letter establishes a novel relationship between a class of recurrent neural networks and certain evolutionary dynamics that emerge in the context of population games. Specifically, it is shown that the output of a recurrent neural network, in the context of classification problems, coincides with the evolution of the population state in a population game. This connection is established via replicator evolutionary dynamics with dynamic payoffs. The connection provides insights into the neural network's behavior from both dynamical systems and game-theoretical perspectives and aligns with recent literature suggesting that the outputs of neural networks may exhibit similarities to the Nash equilibria of suitable games. It also uncovers potential connections between the neural network classification problem and mechanism design. The results are illustrated via different numerical examples.
\end{abstract}


\section{Introduction}
\label{sec:introduction}
\IEEEPARstart{R}{ecurrent} neural networks (RNNs) belong to a class of neural networks (NNs) where the output of a previous step serves as an additional input in the recurrent step. As a result, RNNs are particularly suited for addressing problems with time dependencies. Moreover, they have been recently studied using tools from dynamical systems theory \cite{jafarpour2023efficient}, such as neural ordinary differential equations (Neural ODEs) \cite{Chen2018}, leading to novel applications in control \cite{krstic2023neural} and optimization \cite{liu2021second}. More recently, connections with game theory \cite{marden2018game,brown2017studies} have also been explored, triggering new lines of research for the study and understanding of RNNs \cite{bhatia2014recurrent}.

When discussing the interplay between machine learning, control theory, and game theory, two primary approaches come to the forefront. On one hand, there is a research trend that focuses on leveraging machine learning techniques to tackle complex control and game-related challenges. We shall refer to this approach as \emph{NNs for game theory}. On the other hand, other approaches focus on comprehending the functionality of machine learning techniques through the lens of game theory. This perspective is usually denoted as \emph{game theory for NNs}. Regarding the former approach, numerous works have used machine learning techniques to learn and/or approximate solutions for both control and game theoretical problems \cite{Hu2023}, including optimal control  \cite{Effati2013}, mean-field games \cite{Gomes2023b}, and adversarial behavior in stochastic dynamical systems \cite{BaChoBou2022}.

On the other hand, when studying 
\emph{game theory for NNs}, several recent connections have been uncovered. For instance, in  \cite{Vesseron2021,Ren2021} it was shown that the training of a class of deep NNs is related to the computation of the optimal flow in certain congestion games. Also, in \cite{Tembine2019} and \cite{Jin2020}, Generative Adversarial Networks (GAN) are studied as a min-max game, a connection that is also investigated in \cite{Chivukula2017}, using a two-agent sequential Stackelberg game approach.

This paper narrows its focus to the second approach mentioned earlier, namely game theory for NNs. In particular, we establish a novel relationship between the output of certain RNNs and a class of evolutionary dynamics that emerge in population games \cite{sandholm2010population}: the replicator dynamics. Such dynamics, well-known in the literature of biology and economics \cite{hofbauer1998evolutionary}, have also been used for the synthesis and analysis of feedback controllers and optimization algorithms in different engineering applications \cite{park2019population,quijano2017role,BaTe_2018,arcak2020dissipativity}. Our main contributions are twofold: First, we establish a theoretical link between the output generated by RNN ordinary differential equations (ODEs) in classification problems and the solutions of the replicator dynamics. Specifically, we show that such output follows a replicator system with dynamic payoff functions characterized by the structure of the RNN ODE. We also show that the corresponding population game in the emerging evolutionary dynamics may incorporate payoff dynamics or time-varying payoffs depending on the characteristics of the output layer. This observation highlights potential connections with recent works on population games with dynamic payoffs, and it shows that the task performed by RNNs could similarly be accomplished via evolutionary dynamics and mechanism design. Moreover, it opens the door for the design and analysis of graph-dependent RNNs via tools from graphical population games \cite{pantoja2011population,BaObQu2017,poveda2015shahshahani}. Second, we present experimental validations of the theoretical results in different classification problems, highlighting the connections between the geometry of the replicator dynamics and the properties of the RNN ODE. 

The rest of this letter is organized as follows. Section \ref{sec:background} introduces some preliminaries, Section \ref{sec:bridge} presents the main result, Section \ref{sec:numerical_examples} presents the numerical examples, and Section \ref{sec:conclusions_nn} ends with the conclusions.
\section{Preliminaries}
\label{sec:background}
\subsection{Neural Networks}
\begin{defi}[Set of Layers]
Given $d_{\mathrm{in}},~d_{\mathrm{out}}\in\mathbb{Z}_{\geq1}$, and a map $\sigma:\mathbb{R}\to\mathbb{R}$, a set of \emph{layers} from dimension $d_{\mathrm{in}}$ to $d_{\mathrm{out}}$, with activation function $\sigma$, is given by
\begin{align}\label{eq:layers}
\mathbb{L}^{\sigma}_{[d_{\mathrm{in}},d_{\mathrm{out}}]} &:= \Big\{ \phi : \mathbb{R}^{d_{\mathrm{in}}} \to \mathbb{R}^{d_{\mathrm{out}}}|~\exists~b \in \mathbb{R}^{d_{\mathrm{out}}},~\notag\\
&~~~~~~~\exists~W \in \mathbb{R}^{d_{\mathrm{out}} \times d_{\mathrm{in}}},~\phi_{j} = \sigma\Big(b_j + \sum_{\ell=1}^{d_{\mathrm{in}}} W_{j\ell} x_\ell \Big),~\notag\\
&~~~~~~~\forall~j \in \{1,\dots, d_{\mathrm{out}}\},~\forall~x \in \mathbb{R}^{d_{\mathrm{in}}}\Big\}.
\end{align}
The parameters $b_i$ and $W_{i\ell}$ are known as the biases, and weights, respectively. \QEDB 
\end{defi}

\begin{figure}
    \centering
    \includegraphics[width=0.85\linewidth]{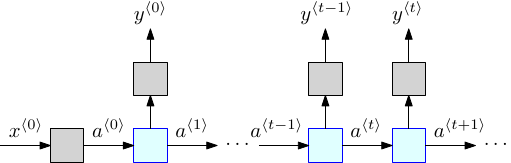}
    \caption{Scheme of a one-to-many single-layer RNN, with $x^{\langle 0\rangle}:=x$.}
    \label{fig:SRNN} 
\end{figure}

Let $\phi^{(i)}\in \mathbb{L}^{\sigma}_{[d_{\mathrm{in}},d_{\mathrm{out}}]}$ for all $i\in\{1,2,\ldots,L\}$, where $L\in\mathbb{Z}_{\geq2}$. Then, the collection of $L$-layers $\{\phi^{(i)}\}_{i=1}^{L}$, has \emph{trainable parameters} given by 
\begin{align}\label{thetaeq}
\theta := \{ W^{(0)}, b^{(0)},W^{(1)}, b^{(1)}, \dots, W^{(L-1)}, b^{(L-1)} \},
\end{align}
where $\theta \in \Theta$, with $\Theta$ being a compact set, and $W^{(t)}$ and $b^{(t)}$ denote the weights and biases of the $t-$th layer, respectively.
\begin{defi}[Set of a Family of Neural Networks]\label{defiNN}
Given $d_{0},~d_{L}\in\mathbb{Z}_{\geq1}$ and $L\in\mathbb{Z}_{\geq2}$, the set of Neural Networks (NNs) with $L$ layers, from dimension $d_{0}$ to $d_{L}$, is given by
\begin{align}\label{eq:family_NN}
\Xi_{[d_0,d_L]} &:= \Big\{ \xi^{\theta}_{[d_0,d_L]} : \mathbb{R}^{d_0} \to \mathbb{R}^{d_L}~|~\exists~\phi^{(0)} \in \mathbb{L}^{\sigma_{\mathrm{i}}}_{[d_0,d_{1}]},\notag\\ 
&~~~~~~~\exists~\phi^{(k)} \in \mathbb{L}^{\sigma_{\mathrm{h}_k}}_{[d_k,d_{k+1}]}, \forall~ k \in \{1,\dots,L-2\},\notag\\
&~~~~~~~\exists~\phi^{(L-1)} \in \mathbb{L}^{\sigma_{\mathrm{o}}}_{[d_{L-1},d_L]}, \notag\\
&~~~~~~~\xi^{\theta}_{[d_0,d_L]} = \phi^{(L-1)} \circ \phi^{(L-2)} \circ \dots \circ \phi^{(0)}
\Big\},
\end{align}
where $\sigma_{\mathrm{i}}$, $\sigma_{\mathrm{h}_k}$, and $\sigma_{\mathrm{o}}$ denote the activation functions of the input layer, the $k-$th hidden layer, and the output layer, respectively. Note that $\theta$ in \eqref{thetaeq} corresponds to the trainable parameters of $\Xi_{[d_0,d_L]}$.
\QEDB
\end{defi}
\begin{remark}
The NNs under study can have any number of hidden layers with different activation functions. Hence, not all the hidden layers need to have the same number of neurons. Later, when working with RNNs, the hidden component can be composed of either a single or multiple heterogeneous layers, provided there is compatibility in the input-output dimensionality, which is essential for an ODE-like representation. \QEDB 
\end{remark}
\begin{remark}
For any $L\in\mathbb{Z}_{\geq2}$, a NN $\xi^{\theta}_{[d_0,d_L]}(\cdot)$ defines the input-to-output mapping $y= \xi^\theta_{[d_0,d_L]}(x)$, which can be written recursively as:
\begin{subequations}
\label{eq:nn_layers_superscript}
\begin{align}
a^{(1)} &= \sigma_{\mathrm{i}} \left(W^{(0)} x + b^{(0)}\right),\\
a^{(k+1)} &= \sigma_{\mathrm{h}_k} \left(W^{(k)} a^{(k)} + b^{(k)}\right), k\in \{1,...,L-2\}, \\
y &= \sigma_{\mathrm{o}} \left(W^{(L-1)} a^{(L-1)} + b^{(L-1)}\right),
\end{align}
\end{subequations}
for all $x \in \mathbb{R}^{d_0}$, and $y \in \mathbb{R}^{d_L}$. We will use a similar structure to describe recurrent neural networks, using \emph{time} instead of \emph{layers} to describe the recursive construction\footnote{In this paper, we use the superscript $(\cdot)$ to denote the hidden layers in NNs, and $\langle \cdot \rangle$ to denote time for the RNNs, i.e., $(\cdot) \ne \langle \cdot \rangle$. See Figures \ref{fig:SRNN}-\ref{fig:DRNN}.}. \QEDB 
\end{remark}
\subsection{Recurrent Neural Networks and ODE Approximation}
\begin{figure}
    \centering
    \includegraphics[width=0.8\linewidth]{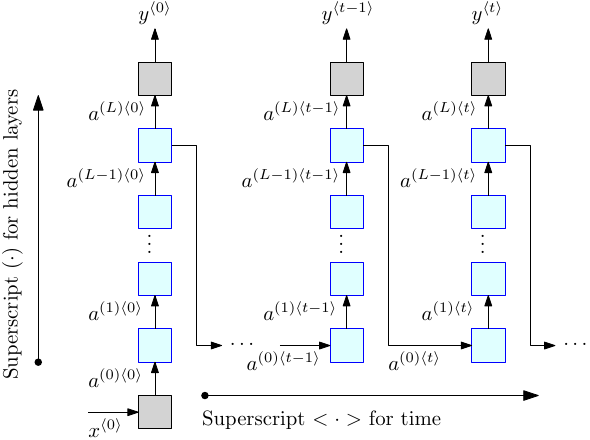}
    \caption{Scheme of a one-to-many multiple-layer RNN, with $x^{\langle0\rangle}:=x$.}
    \label{fig:DRNN}
\end{figure}

There are multiple classes of RNNs, primarily distinguished by the input/output information they receive/return, as well as the architecture of their hidden layers. We focus on two particular classes of RNNs: (i) Single-Layer (shallow) one-to-many, shown in Figure \ref{fig:SRNN}, and (ii) Multiple-Layer (deep) one-to-many, shown in Figure \ref{fig:DRNN}.  
\begin{defi}
An RNN is said to be of Class I if its input-to-output relation satisfies the following time-recursive relations for all $x\in\mathbb{R}^M$:
\begin{subequations}\label{eq:RNN_1}
\begin{align}
a^{\langle 0 \rangle} &= \sigma_{\mathrm{i}} \left(W_x x^{\langle 0 \rangle} + b_x\right),~~x^{\langle0\rangle}=x,\label{eq:RNN_1a}\\
a^{\langle t+1 \rangle} &= \sigma_{\mathrm{h}} \left(W_a a^{\langle t \rangle} + b_a\right),~~t = 0,1,\dots \label{eq:RNN_1b}\\
y^{\langle t \rangle } &= \sigma_{\mathrm{o}} \left(W_y a^{\langle t \rangle} + b_y\right),~~t = 0,1,\dots \label{eq:RNN_1c}
\end{align}   
\end{subequations}
where $a^{\langle t\rangle} \in \mathbb{R}^{d_0}$, and $y^{\langle t\rangle} \in \mathbb{R}^N$, for all $t\in\mathbb{Z}_{\geq0}$. \QEDB 
\end{defi}
\begin{defi}
An RNN is said to be of Class II if its input-to-output relation satisfies the following time-recursive relations for all $x\in\mathbb{R}^M$:
\begin{subequations}
\label{eq:RNN_2}
 \begin{align}
a^{\langle 0 \rangle} &= \sigma_{\mathrm{i}} \left(W_x x^{\langle 0 \rangle} + b_x\right),~~x^{\langle0\rangle}=x \label{eq:RNN_2a}\\
a^{\langle t+1 \rangle} &= \xi^{\theta}_{[d_0,d_{L}]}\left(a^{\langle t \rangle}\right),~~~~~~~t = 0,1,\dots  \label{eq:RNN_2b}\\
y^{\langle t \rangle} &= \sigma_{\mathrm{o}} \left(W_y a^{\langle t \rangle} + b_y\right),~~t = 0,1,\dots  \label{eq:RNN_2c}
\end{align}   
\end{subequations}
where $d_0=d_L$. \QEDB 
\end{defi}

 Note that \eqref{eq:RNN_1}-\eqref{eq:RNN_2} differ from \eqref{eq:nn_layers_superscript} in the superscript, see also Figures \ref{fig:SRNN} and \ref{fig:DRNN}. We consider RNNs that satisfy the following: 
\begin{assump}\label{regularityassumption}
The maps $\sigma_i(\cdot),\sigma_h(\cdot),\sigma_o(\cdot)$ and $\xi^{\theta}_{[d_{0},d_L]}(\cdot)$ are Lipschitz continuous. \QEDB 
\end{assump}

The recursive updates in \eqref{eq:RNN_1} or \eqref{eq:RNN_2} can be interpreted as an Euler discretization of a suitable ODE \cite{Chen2018}. To illustrate this idea, we focus on Class II RNNs, but the approach is also easily applicable to RNNs of Class I. Let $\mathcal{T} := \{0,\dots,N_T\}$ denote a discrete-time window used to implement the RNN \eqref{eq:RNN_2}. Using a small ``step size'' parameter $\tau=T/N_T$ for a continuous-time interval $[0,T]$, where $T>0$, let
\begin{align}\label{discretization}
\tilde{\xi}^{\theta}(a{( t )})&:= \frac{a^{\langle t+1 \rangle} -  a^{\langle t \rangle}}{\tau} =  \frac{\xi^{\theta}_{[d_0,d_{L}]}(a^{\langle t \rangle}) -  a^{\langle t \rangle} }{\tau},
\end{align}
where $a^{\langle t+1 \rangle}$ is defined in \eqref{eq:RNN_2}. Thus, when $N_T$ is sufficiently large, i.e., $\tau$ is sufficiently small, we can approximate the recursive dynamics of the hidden layers within the discrete-time window $\mathcal{T}$ via the following Lipschitz ODE 
\begin{subequations}\label{eq:RNN_shadow}
\begin{align}
\dot{a}{( t )} &= \tilde{\xi}^{\theta}(a{( t )}),~~~~~~a{( 0 )}= \phi^{\mathrm{in}}(x),\\
y{( t )} &= \phi^{\mathrm{out}}(a(t)),
\end{align}
\end{subequations}
where the input and output layers are:
\begin{subequations}
\begin{align}
\phi^{\mathrm{in}}(x) &=\sigma_{\mathrm{i}} \left(W_x x + b_x\right),~~~~~\forall~x\in\mathbb{R}^M,\\
\phi^{\mathrm{out}}(a(t)) &= \sigma_{\mathrm{o}} \left(W_y a{( t )} + b_y\right),~~\forall~t\in[0,T],
\end{align}
\end{subequations}
with $\phi^{\mathrm{in}} \in \mathbb{L}^{\sigma_i}_{[M,d_0]}$ and $\phi^{\mathrm{out}} \in \mathbb{L}^{\sigma_o}_{[d_L,N]}$.
We refer to \eqref{eq:RNN_shadow} as the RNN ODE, and its trainable parameters are given by
\begin{align}
\label{eq:RNN_parameters}
\hat{\theta} := (W_x,b_x,W_y,b_y,\theta),
\end{align}
where $\theta$ is given by \eqref{thetaeq}. Under Assumption \ref{regularityassumption}, the Euler discretization \eqref{discretization} provides a suitable approximation for \eqref{eq:RNN_shadow} on compact sets, provided $\tau$ is sufficiently small \cite[Thm. 5.2]{sanfelice2010dynamical}.
\subsection{The Classification Problem via RNN}
The classification problem consists of assigning to an unknown input datum, denoted by $x$, a label from a finite set of possibilities $\mathcal{L}=\{1,\dots,N\}$. The output of the neural network, denoted by $y$, is interpreted as a prediction of the true label corresponding to $x$. To train the neural network, true information is provided in a data set $\mathcal{D} = \{1,\dots,D\}$. Using these true labels, for each $x^j$ denoting the $j-$th input training datum for which the true label is known, we let
\begin{align*}
l_i^j = \begin{cases}
1,~\text{if}~x^j~\text{is of label}~i,\\
0,~\text{otherwise}.
\end{cases}
\end{align*}
for all $j\in \mathcal{D},~i \in \mathcal{L}$, and we let
$l^j:= [l^j_1 \quad \dots \quad l^j_N]^{\top}$. Based on this, we consider the following loss function: 
\begin{align*}
L(\hat{\theta}) = - \sum_{k\in \mathcal{T}} \sum_{j \in \mathcal{D}} \sum_{i \in \mathcal{L}} l^j_i \log(y^{j\langle k \rangle}_{i}).
\end{align*}
where $y^{j\langle k \rangle}_{i}$ is given by \eqref{eq:RNN_2c}, i.e., it is the output of the RNN at time $k$, of the $j-$th training datum corresponding to the $i-$th label. Thus, by using the ``softmax'' function, the RNN training problem that we study can be formally stated as:
\begin{align*}
\min_{\hat{\theta}}~~L(\hat{\theta}),~~\text{s.t.~} 
\left\{\begin{array}{l}
 a^{j\langle k+1 \rangle}= a^{j\langle k \rangle} + \tau \cdot \tilde{\xi}^{\theta}(a^{j\langle k \rangle}),\\
~~~y^{j\langle k \rangle}= \mathrm{softmax}(a^{j\langle k \rangle}),\\
~~~a^{j\langle 0 \rangle}= \phi^{\mathrm{in}}(x^{j\langle 0 \rangle}),\\
~~~x^{j\langle 0 \rangle}= x^{j},~~\forall~ j \in \mathcal{D},~k \in \mathcal{T}.
\end{array}\right.
\end{align*}
From the dynamical systems point of view, each data point creates a trajectory of length $T$ for the RNN ODE. Therefore, by solving the above optimization problem the goal is to make the trajectory converge to the corresponding true label of the data point.
\subsection{Population games and Replicator Systems}
\label{sec:population_dynamics}
To bridge the earlier discussion with the replicator system, consider a large population of decision-makers who can select a strategy from a discrete set $\mathcal{L}:= \{1,\dots,N\}$. Let $p_i$ be the proportion of entities that choose the strategy $i$. Then, the vector $p=(p_1,p_2,\ldots,p_N)$, usually referred to as the ``state'' of the population \cite{sandholm2010population}, satisfies
\begin{equation}
p\in \Delta := \{ p \in \mathbb{R}^N_+: \mathbf{1}_N^\top p = 1 \},
\end{equation}
where $\Delta$ is the standard unitary simplex. Let $f:\Delta\to\mathbb{R}^N$ be a vector of \emph{payoffs} that describes the desirability of each of the strategies in the set $\mathcal{L}$. Thus, the $i-$th component of $f$, denoted $f_i$, describes the ``payoff'' related to the $i$th strategy. A typical model in biology and economics that captures the evolution over time of the population state $p$ in the above setting is described by the replicator dynamics \cite{sandholm2010population}: 
\begin{figure}[t!]
\centering
\includegraphics[width=0.23\textwidth]{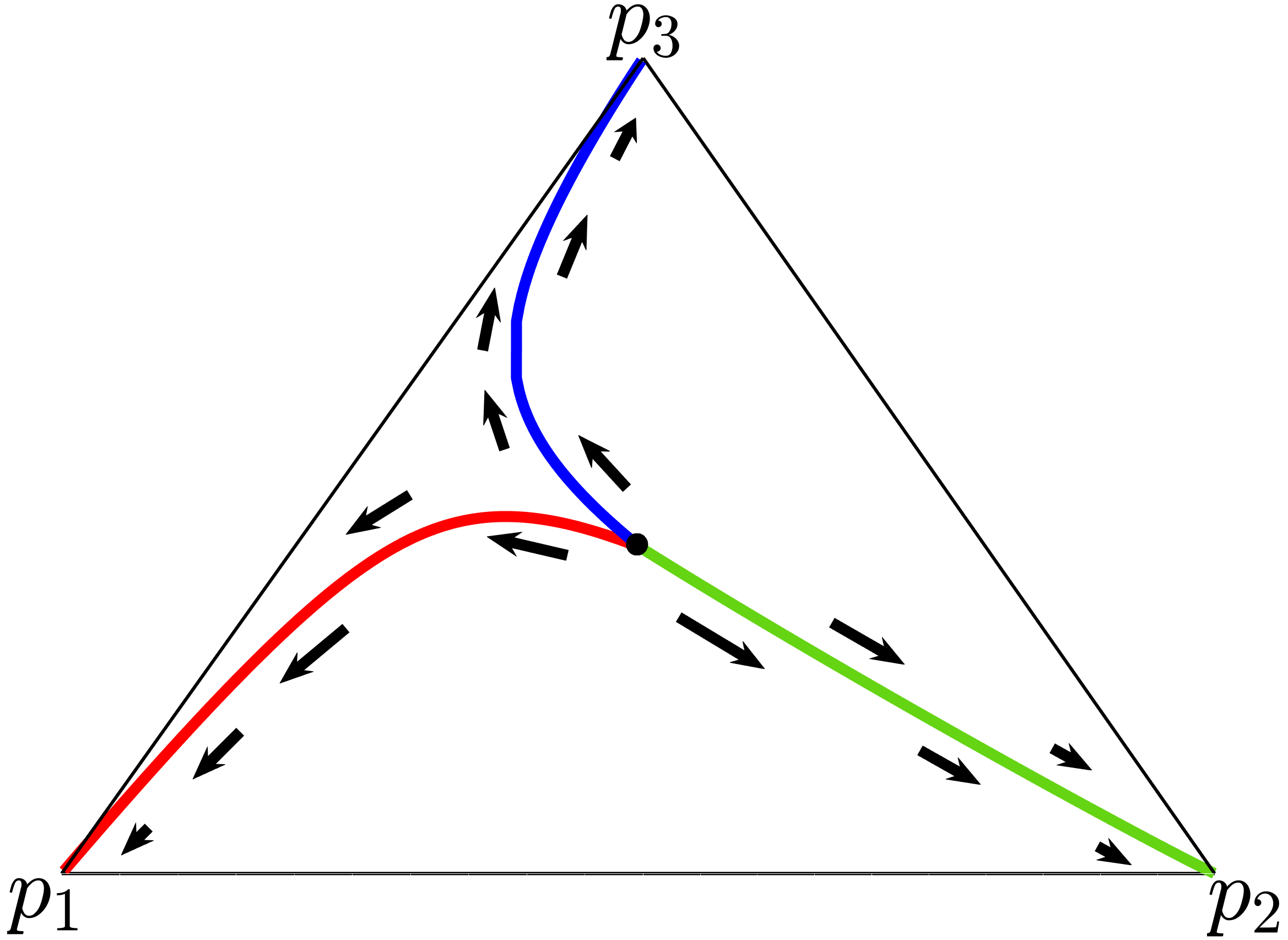}
\includegraphics[width=0.23\textwidth]{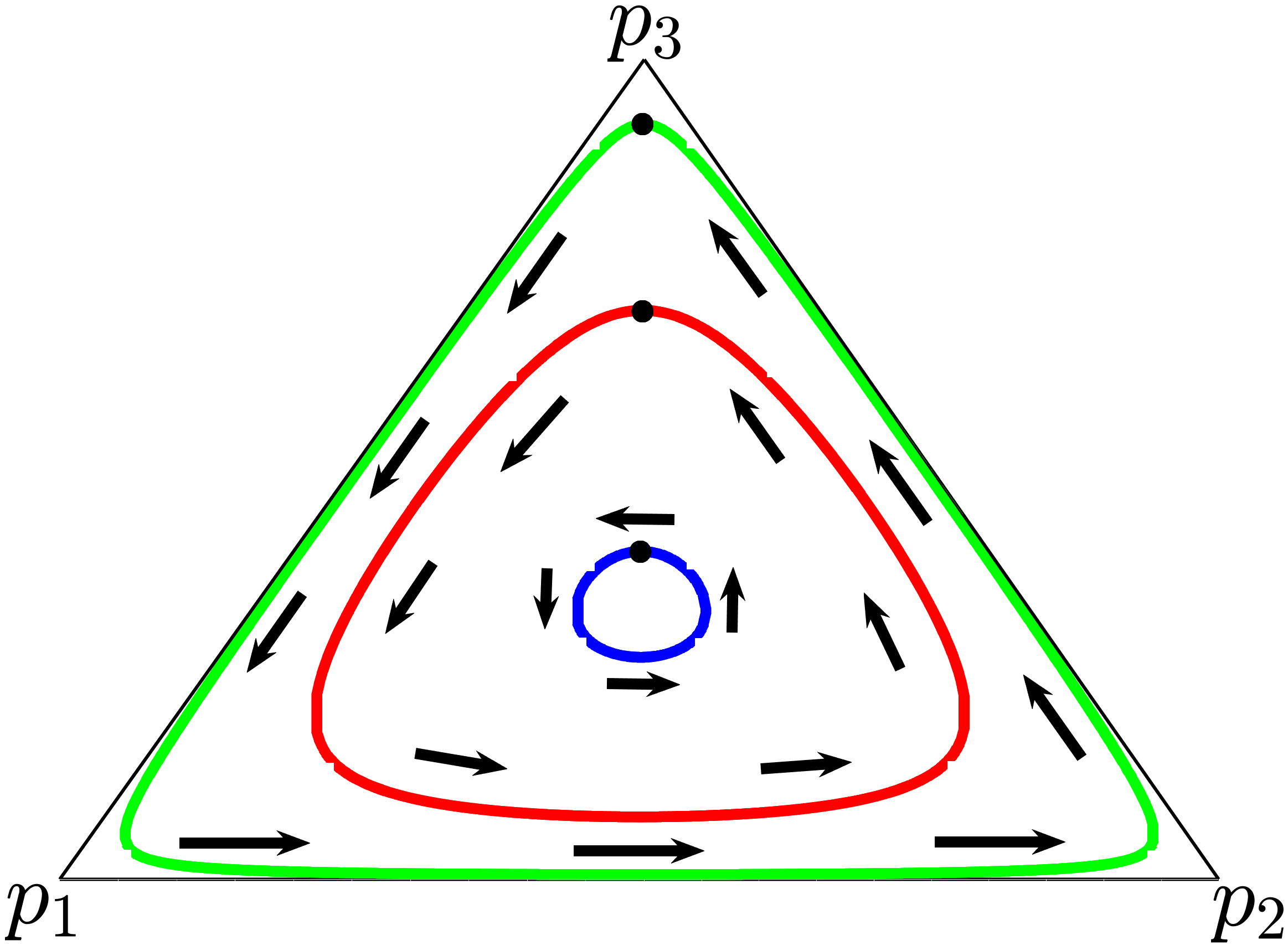}
  \caption{Different trajectories, converging (left) and oscillating (right), of the replicator dynamics \eqref{replicatorequations} under different choices of $f$, evolving on the simplex with $N=3$. The pure strategies $\mathcal{L}$ correspond to the corners of the simplex. The black dot indicates the initial conditions.}
  \label{fig:scheme1}
\end{figure}
\begin{align}\label{replicatorequations}
	\dot{p}{(t)} &= \mathrm{diag}(p{( t )}) \Big( f(p( t )) - \mathbf{1}_N  f(p( t ))^\top p(t) \Big),
\end{align}
with $p{(0)} \in \Delta$. Replicator systems of the form \eqref{replicatorequations} have been studied in biology \cite{hofbauer1998evolutionary}, economics \cite{sandholm2010population}, and control engineering systems \cite{BaObQu2017}. Moreover, when $f$ is the gradient of a potential field, system \eqref{replicatorequations} describes a gradient flow in a non-Euclidean metric. These connections have linked the replicator dynamics with adaptive systems able to achieve self-optimizing behaviors \cite{BaObQu2017,poveda2015shahshahani}. While the scope of this paper does not include a detailed review of the replicator system, it is worth mentioning some interesting properties of \eqref{replicatorequations}. For instance, the replicator dynamics satisfies a property called \textit{Nash stationarity}, meaning that its rest points correspond to Nash equilibria (i.e., every strategy in use earns the maximal payoff) of a game defined by the vector of payoffs $f$ in the simplex $\Delta$, \cite[pp. 24]{sandholm2010population}. Additionally, the dynamics render the set $\Delta$ forward invariant, meaning that if the initial population state is a distribution of agents along the strategies, then a well-defined distribution will be maintained for all time $t\geq0$. Figure \ref{fig:scheme1} presents different sample trajectories of \eqref{replicatorequations} in games that generate converging (left plot) and oscillating (right plot) behaviors. The trajectories are shown evolving over the simplex $\Delta$ when $N=3$. Similar plots will be shown to emerge in Section \ref{sec:numerical_examples} in the context of RNN ODEs.

\section{Main Result}
\label{sec:bridge}
The following theorem is the main result of this letter. Numerical examples are presented in Section \ref{sec:numerical_examples}.
\begin{theorem}\label{propos:replicator1}
Consider a trained RNN $\xi^{\theta^*}_{[d_0,d_{L}]}:\mathbb{R}^{d_0} \to \mathbb{R}^{d_L}$ that satisfies Assumption \ref{regularityassumption}, with parameters $$\theta^* := \{ W^{(0)*}, b^{(0)*},W^{(1)*}, b^{(1)*}, \dots, W^{(L-1)*}, b^{(L-1)*} \},$$ and trained input/output layers given by $\phi^{\mathrm{in}} \in \mathbb{L}^{\sigma}_{[M,d_0]}$ and $\phi^{\mathrm{out}} \in \mathbb{L}^{\sigma}_{[d_L,N]}$, respectively. Then, the following holds:
\begin{enumerate}
\item The composition $y=\phi^{\mathrm{out}} \circ \xi^{\theta^*}_{[d_0,d_{L}]} \circ \phi^{\mathrm{in}}(x)$ returns a classification $y \in \Delta \subset [0,1]^N $ of the input $x \in \mathbb{R}^{M}$ into a label from the set $\mathcal{L} = \{1,\dots,N\}$.
\item  The function $t\mapsto y(t)$ satisfies the following replicator system for all $t\in[0,T]$:
\begin{subequations}\label{firstRDs}
\begin{align}
\dot{y}{( t )} &= \mathrm{diag}(y{( t )}) \left( f(a(t)) - \mathbf{1}_N y(t)^\top f(a(t)) \right),\\
y{(0)} &= \phi^{\mathrm{out}} \circ \phi^{\mathrm{in}}(x),~y{(0)} \in \Delta,
\end{align}
\end{subequations}
with dynamic payoffs $f(a(t)) = W_y \tilde{\xi}^{\theta}(a{( t )})$, where
\begin{subequations}\label{adynamics1}
\begin{align}
\dot{a}{( t )}&= \tilde{\xi}^{\theta}(a{( t )})= \frac{1}{\tau} \left( \xi^{\theta}_{[d_0,d_{L}]}(a^{\langle t \rangle}) -  a^{\langle t \rangle} \right),\\
a{(0)}&= \phi^{\mathrm{in}}(x).
\end{align}
\end{subequations}
\item If the output layer is invertible, then the function $t\mapsto y(t)$ satisfies the following for all $t\in[0,T]$:
\begin{subequations}
\begin{align}
\dot{y}{( t )} &= \mathrm{diag}(y{( t )}) \left( f(y( t ),t) - \mathbf{1}_N  f(y( t ),t)^\top y(t) \right),\\
y{(0)} &= \phi^{\mathrm{out}} \circ \phi^{\mathrm{in}}(x)\in \Delta,
\end{align}
\end{subequations}
with payoffs $f$ given by \eqref{timevaryingpayoffs}. \QEDB 
\end{enumerate}
\end{theorem}

\textit{Proof:} The fact that $y$ returns a classification $y\in\Delta$ is directly induced by the softmax function used in the output layer. On the other hand, to establish item (2), let:
\begin{equation}\label{zform}
z(t):= W_y a{(t)} + b_y,~~~\forall~t\geq0,
\end{equation}
where $W_y,b_y$ are the weights of the output layer of the RNN. It follows that the $i-$th component of the output $y$ satisfies 
\begin{equation}\label{outputrepresentation}
y_i(t) = \phi_i^{\mathrm{out}}(a(t))=\mathrm{softmax} \left(z_i(t)\right)=\frac{e^{z_{i}(t)}}{\sum_{j}~e^{z_{j}(t)}},
\end{equation}
for all $t\geq0$. 
Computing the time-derivative of $y_i$, using $\omega_i(t)= e^{z_i(t)}$ for all $i \in \{1,\dots,N\}$, we obtain
\begin{align*}
&\dot{y}_i(t) = \frac{\omega_i(t) \dot{z}_i(t) \sum\limits_{j =1}^N \omega_j(t) - \sum\limits_{j =1}^N \omega_j(t) \dot{z}_j(t) \omega_i(t) }{\left(\sum\limits_{k =1}^N \omega_k(t) \right)^2},\\
&= \frac{\omega_i(t)}{\sum\limits_{k =1}^N \omega_k(t)} 
\left( \dot{z}_i(t) \frac{\sum\limits_{j =1}^N \omega_j(t)}{\sum\limits_{k =1}^N \omega_k(t)}  - \sum_{j =1}^N \left[ \frac{\omega_j(t)}{\sum\limits_{k =1}^N \omega_z(t)} \dot{z}_j(t) \right]\right).
\end{align*}
Using $y_i(t)=\frac{\omega_i(t)}{\sum_{k=1}^N\omega_k(t)}$ and $\frac{\sum_{j=1}^N\omega_{j}(t)}{\sum_{k=1}^N\omega_{k}(t)}=1$, we obtain
\begin{equation}
\dot{y}_i= y_i(t) \left( \dot{z}_i(t) - \sum\limits_{j =1}^N y_j(t) \dot{z}_j(t) \right).
\end{equation}
\begin{figure}[t!]
\centering
\includegraphics[width=0.4\textwidth]{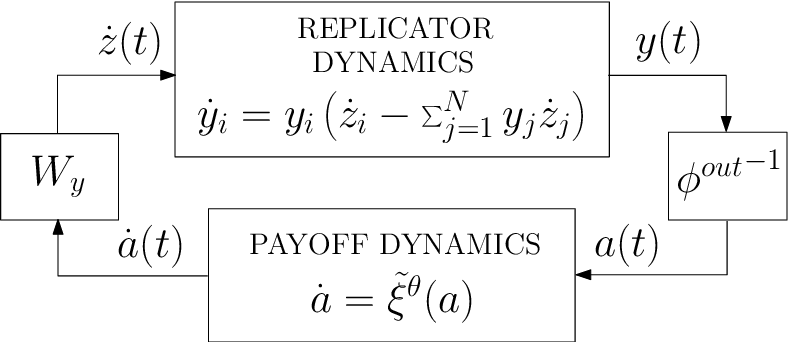} 
  \caption{Feedback interconnection modeling the RNN ODE for output classification problems.}
  \label{fig:scheme}
\end{figure}
It follows that the output vector satisfies
\begin{align*}
\dot{y}{( t )} &= \mathrm{diag}(y{( t )}) \left( \dot{z}{( t )} - \mathbf{1}_N y(t)^\top \dot{z}{( t )} \right).
\end{align*}
Using $\dot{z}(t)=W_y \dot{a}{( t )}$, and the RNN ODE \eqref{eq:RNN_shadow}, we finally obtain the complete dynamics of the system:
\begin{align*}
    \dot{y}{( t )} &= \mathrm{diag}(y{( t )}) \left( W_y \tilde{\xi}^{\theta}(a{( t )}) - \mathbf{1}_N y(t)^\top W_y \tilde{\xi}^{\theta}(a{( t )}) \right),\\
    \dot{a}{( t )} &= \tilde{\xi}^{\theta}(a{( t )}),
\end{align*}
with initial conditions 
\begin{equation}
y{(0)}= \phi^{\mathrm{out}} \circ \phi^{\mathrm{in}}(x)\in \Delta,~~~a(0)= \phi^{\mathrm{in}}(x).
\end{equation}
To establish item (3), suppose now that the output layer is invertible, i.e., the inverse of the matrix $W_y$ exists. Using again $\omega_i(t)= e^{z_i(t)}$, for all $i \in \{1,\dots,N\}$, we have
\begin{align*}
    \ln(y_i(t)) &= z_i(t) -\ln\left(\sum_{k =1}^N \omega_k(t)\right)=z_i(t)-C(t),
\end{align*}
which implies that $z_i(t)=  \ln(y_i(t)) + {C(t)}$. Using \eqref{zform} and solving for $a$, we have
\begin{align*}
     a{( t )}  &= W_y^{-1} \ln(y(t)) + W_y^{-1} \left( \mathbf{1}_N{C(t)}- b_y \right),
\end{align*}
leading to the following replicator system:
\begin{align*}
    \dot{y}{( t )} &= \mathrm{diag}(y{( t )}) \left( f(y(t),t) - \mathbf{1}_N y(t)^\top f(y(t),t) \right),
\end{align*}
with vector of payoff functions
\begin{align}\label{timevaryingpayoffs}
    f(y(t),t) = W_y \tilde{\xi}^{\theta}\left( W_y^{-1} \left(\ln(y(t)) + \mathbf{1}_N{C(t)}- b_y \right) \right)
\end{align}
and initialization $y(0)= \phi^{\mathrm{out}} \circ \phi^{\mathrm{in}}(x)\in \Delta$. \hfill $\blacksquare$
\begin{remark}
If there exists a closed-form solution $t\mapsto a(t)$ to the RNN ODE $a(t)=\phi^{\text{in}}(x)+\int_{0}^{t}\tilde{\xi}^{\theta}(a(\tau))d\tau$,~$t\geq 0$, then substituting in \eqref{firstRDs} leads to the ``standard'' replicator system \eqref{replicatorequations} with time-varying payoffs $f(t)=W_y \tilde{\xi}^{\theta}(a(t))$. \QEDB 
\end{remark}
\begin{figure*}[t!]
	\begin{center}
 \resizebox{\textwidth}{!}{
		\begin{tabular}{cc|cc|cc}
			\multicolumn{2}{c}{\textbf{Experiment 1}} & \multicolumn{2}{c}{\textbf{Experiment 2}} & \multicolumn{2}{c}{\textbf{Experiment 3}} \\ 
			\includegraphics[width=0.2\textwidth]{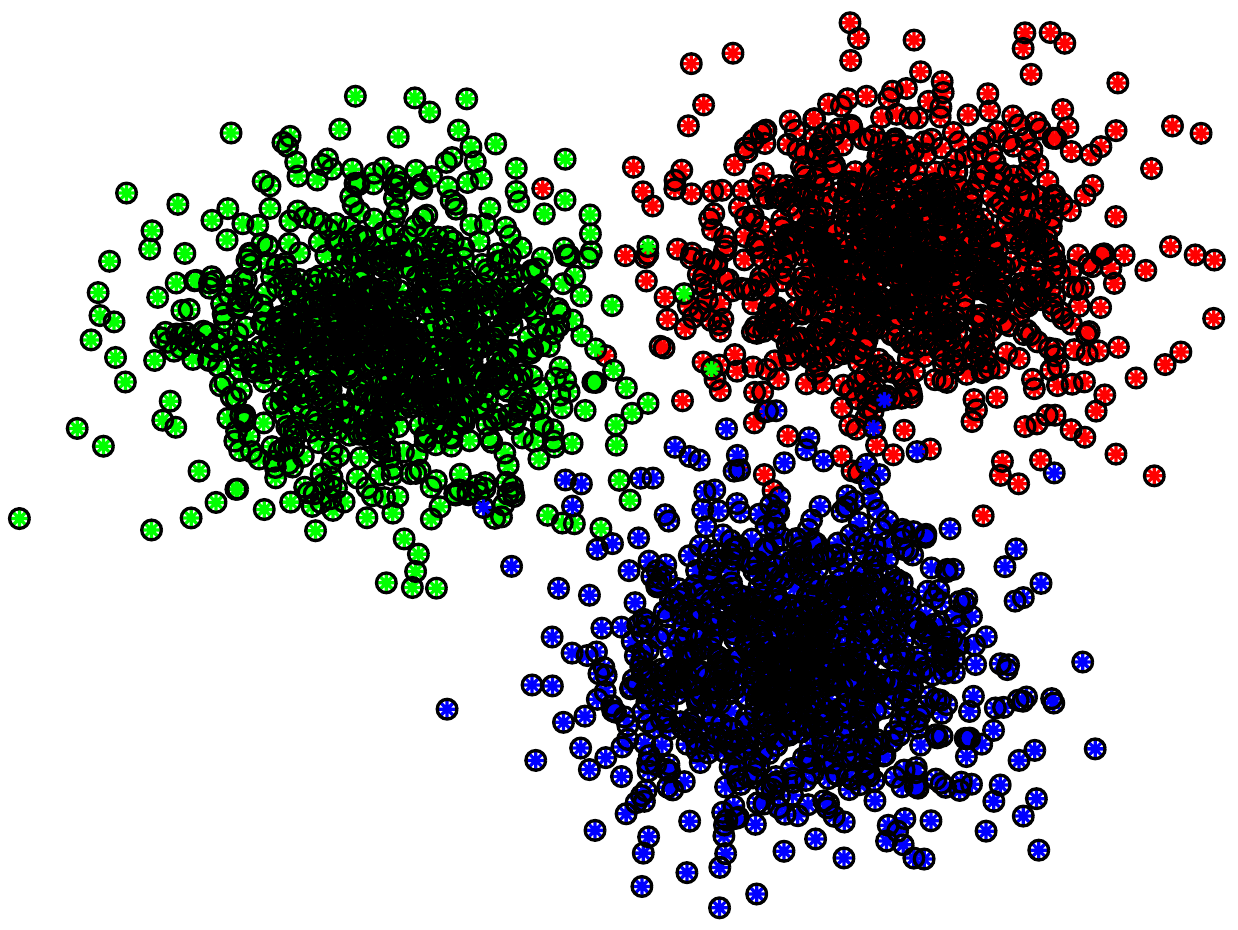} & \includegraphics[width=0.07\textwidth]{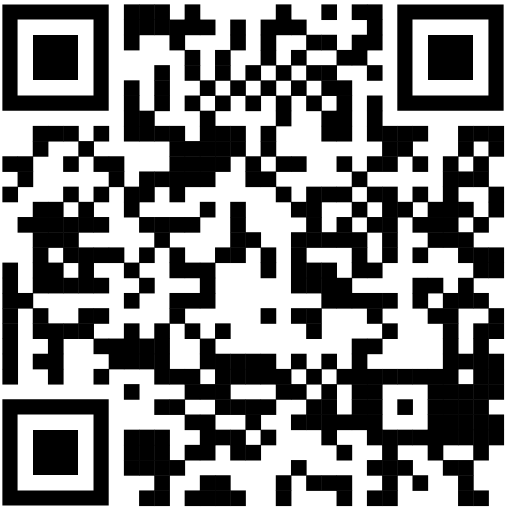} & \includegraphics[width=0.2\textwidth]{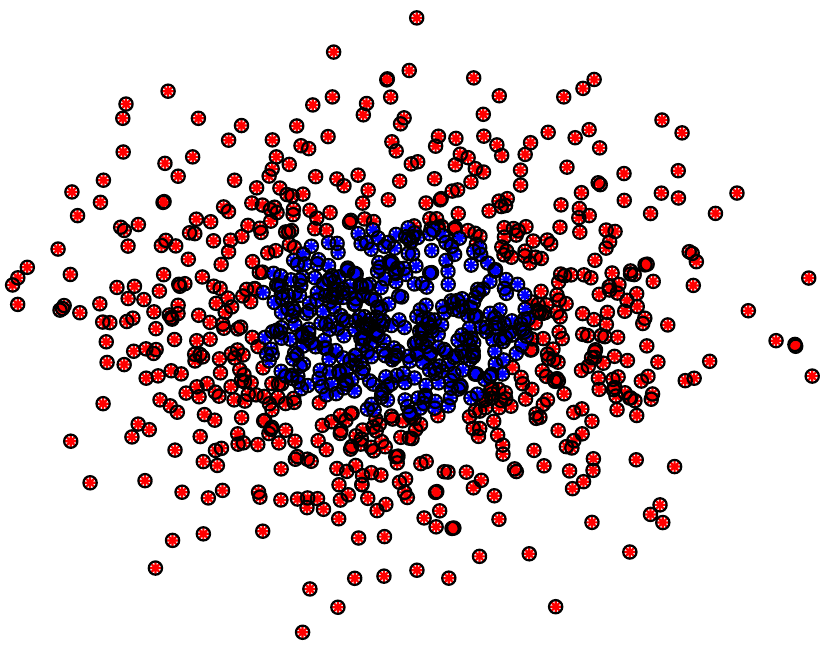}  & \includegraphics[width=0.07\textwidth]{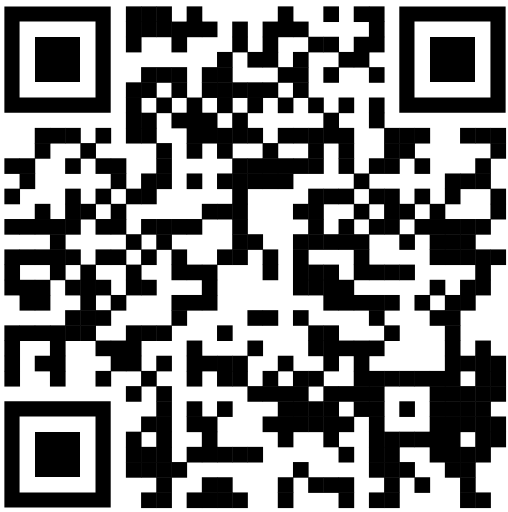} & \includegraphics[width=0.2\textwidth]{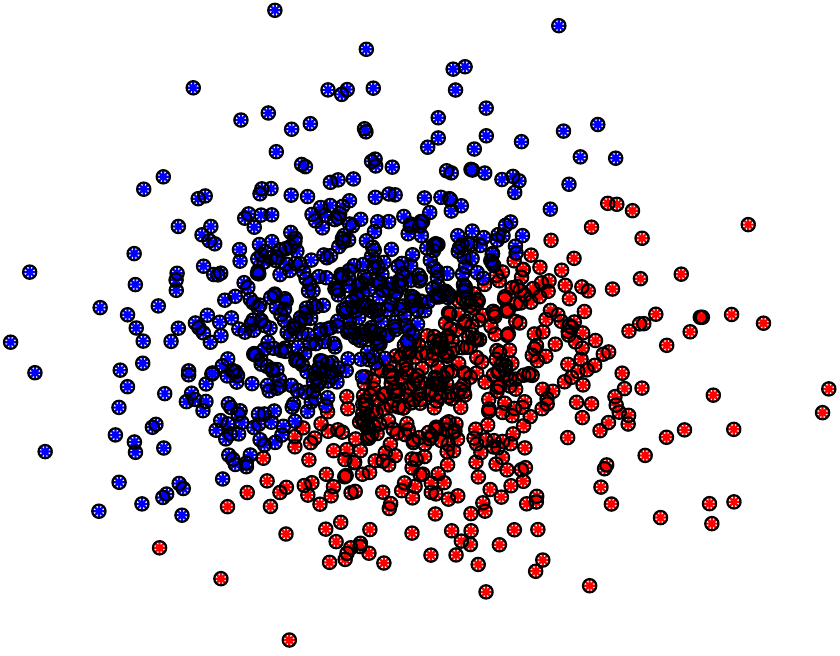} & \includegraphics[width=0.07\textwidth]{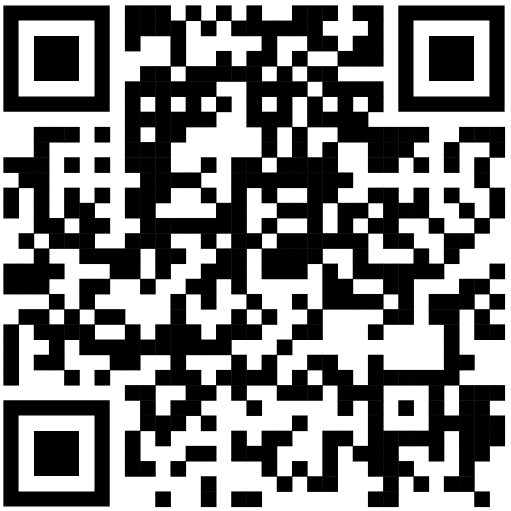} \\
			 $N=3$, $D=3000$ & \textbf{Animated} & $N=2$, $D=1000$ &  \textbf{Animated} & $N=2$, $D=1000$ & \textbf{Animated}\\
			 $x^j \in \mathbb{R}^{2},\forall j \in \mathcal{D}$ & \textbf{Simulation} & $x^j \in \mathbb{R}^{2},\forall j \in \mathcal{D}$ & \textbf{Simulation} &
			$x^j \in \mathbb{R}^{2},\forall j \in \mathcal{D}$ & \textbf{Simulation}
		\end{tabular}
  }
	\end{center}
 \caption{Training datasets for experiments 1, 2, and 3. The QR code can be used to access animated simulations of the system.}
 \label{fig:experiment_data}
\end{figure*}
\begin{remark}
Figure \ref{fig:scheme} presents a block diagram illustrating the interconnection between the replicator dynamics and the payoff dynamics $\dot{a}$. Such feedback systems have been recently studied in the population games literature \cite{arcak2020dissipativity,quijano2017role,poveda2015shahshahani}. Its emergence in the context of RNNs for classification problems suggests potential connections with other neural ODEs utilizing the Boltzmann distribution at the output layer, as well as decentralized dynamical systems defined over networks \cite{pantoja2011population}. When the output layer is non-invertible, system \eqref{firstRDs}-\eqref{adynamics1} can be viewed as a cascade system. The emerging behavior of \eqref{firstRDs} can be ``shaped" through a suitable design of the payoff dynamics \eqref{adynamics1}. Further studies on such shaping mechanisms represent interesting future research directions.  \QEDB 
\end{remark}
\begin{figure*}[h!]
\begin{center}
	\begin{tabular}{cc}
	\includegraphics[width=0.48\textwidth]{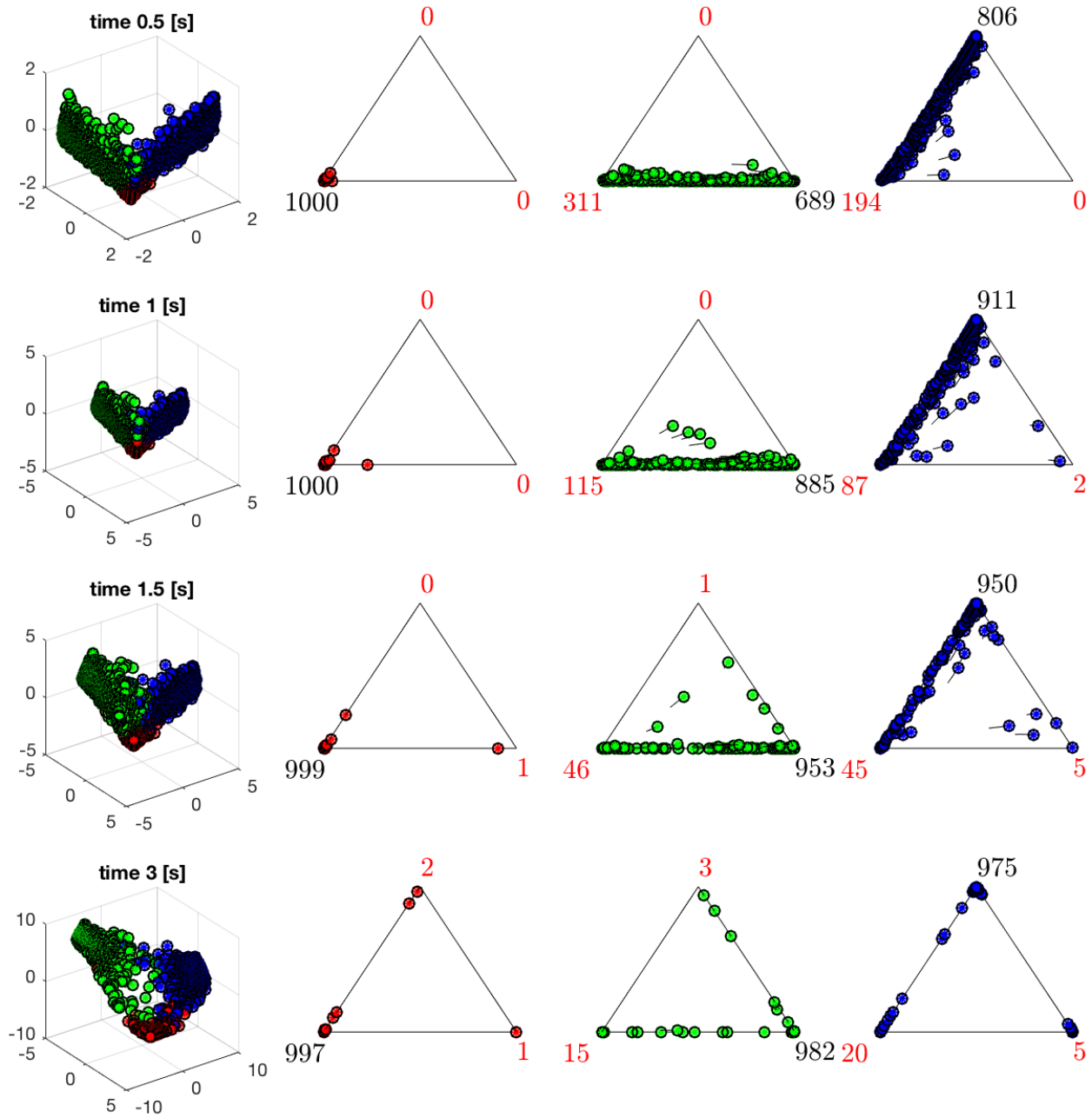} & \includegraphics[width=0.48\textwidth]{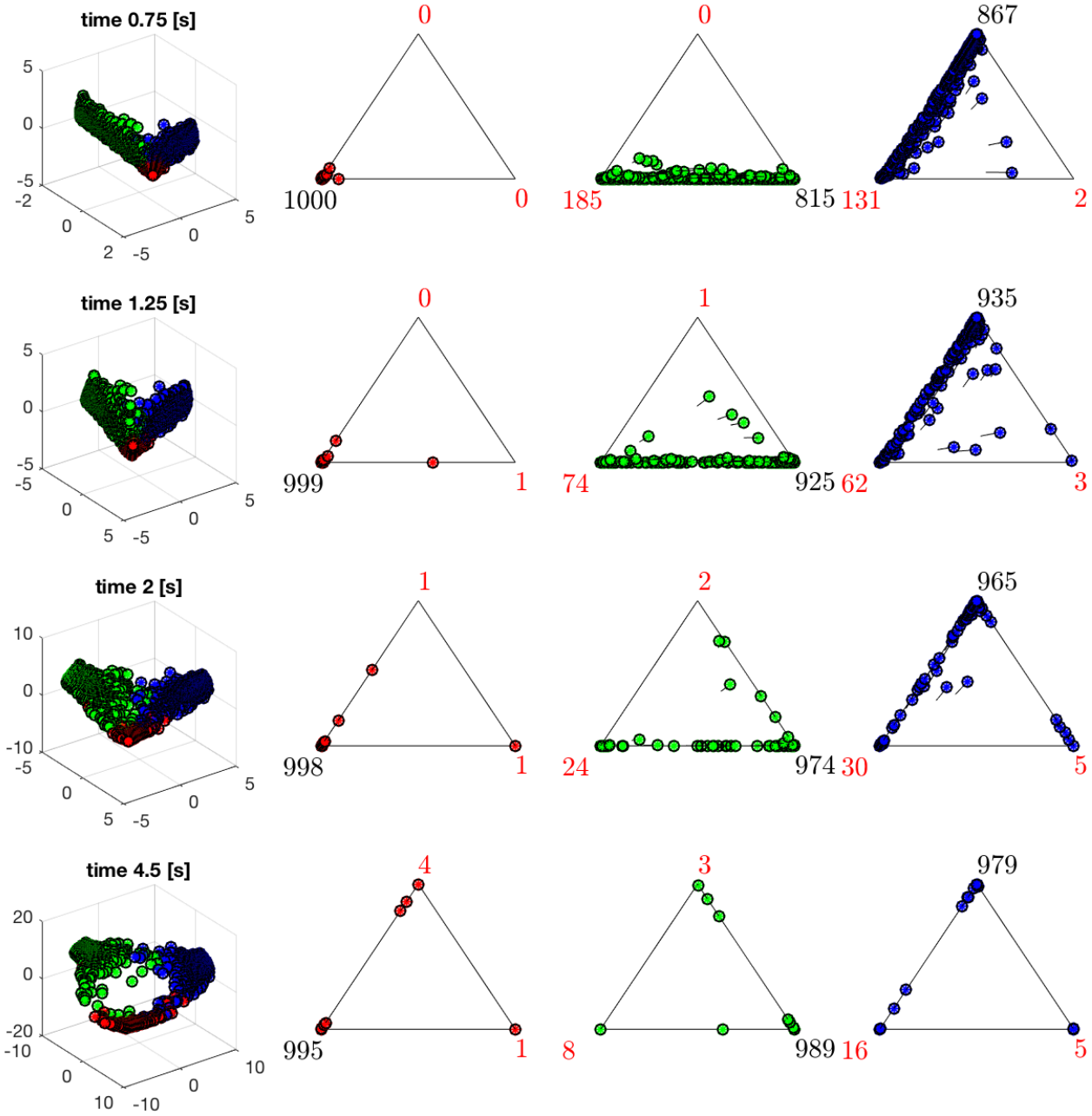}
	\end{tabular}
 \end{center}
 \caption{Results for experiment 1. Temporal evolution of the ODE RNN $\tilde{\xi}^\theta(a(t))$ and the corresponding replicator dynamics for classification in $\Delta$.}
 \label{fig:experiment_1}
\end{figure*}

\section{Numerical Experiments}
\label{sec:numerical_examples}
In this section, we present two numerical examples to illustrate our results.
\subsection{Classification of Data Using 3 Labels}
To illustrate Theorem 1, we first consider three simple low-dimensional classification experiments. 
\begin{figure*}[t!]
\begin{center}
	\begin{tabular}{cc}		\includegraphics[width=0.48\textwidth]{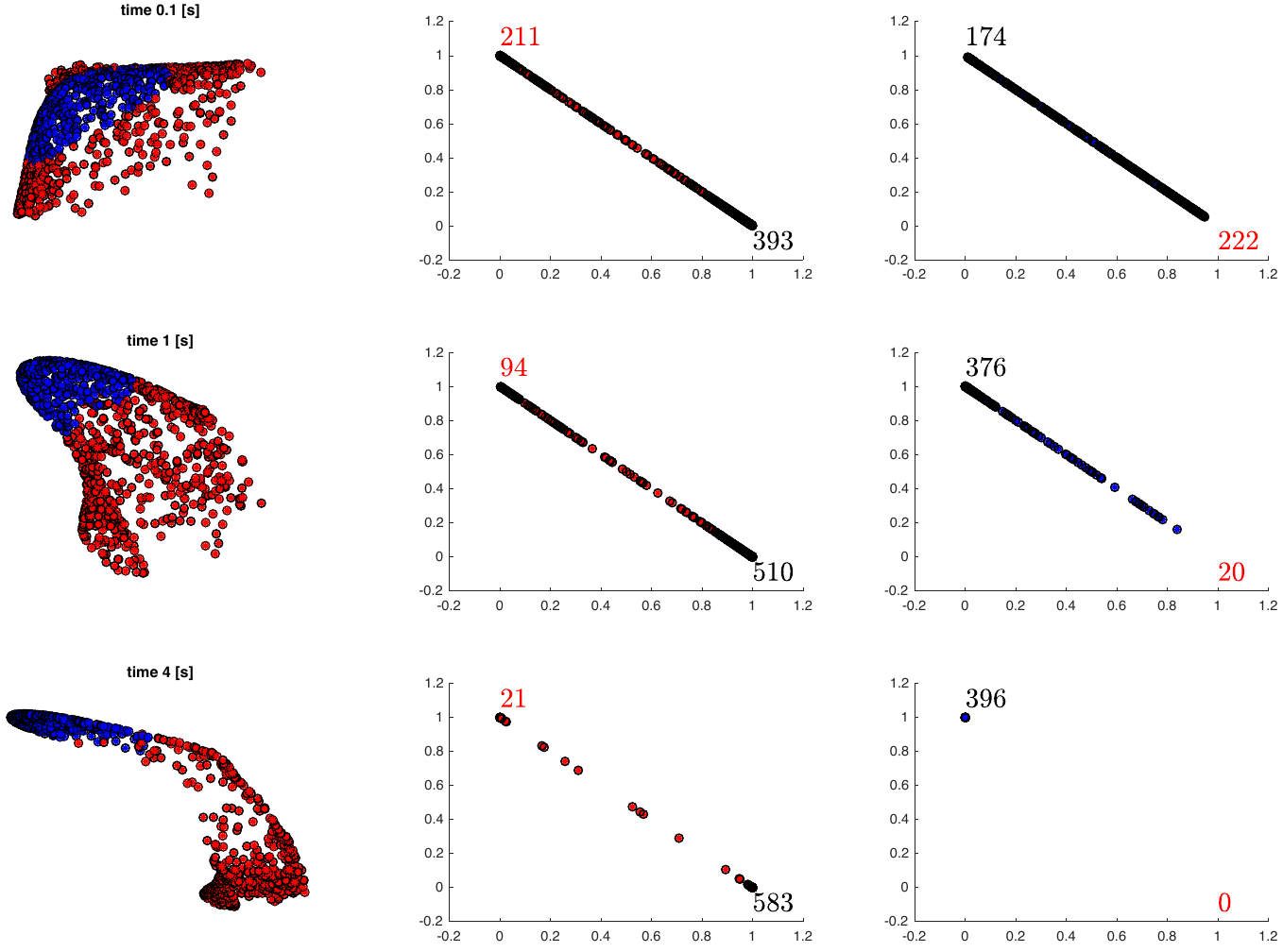} & \includegraphics[width=0.48\textwidth]{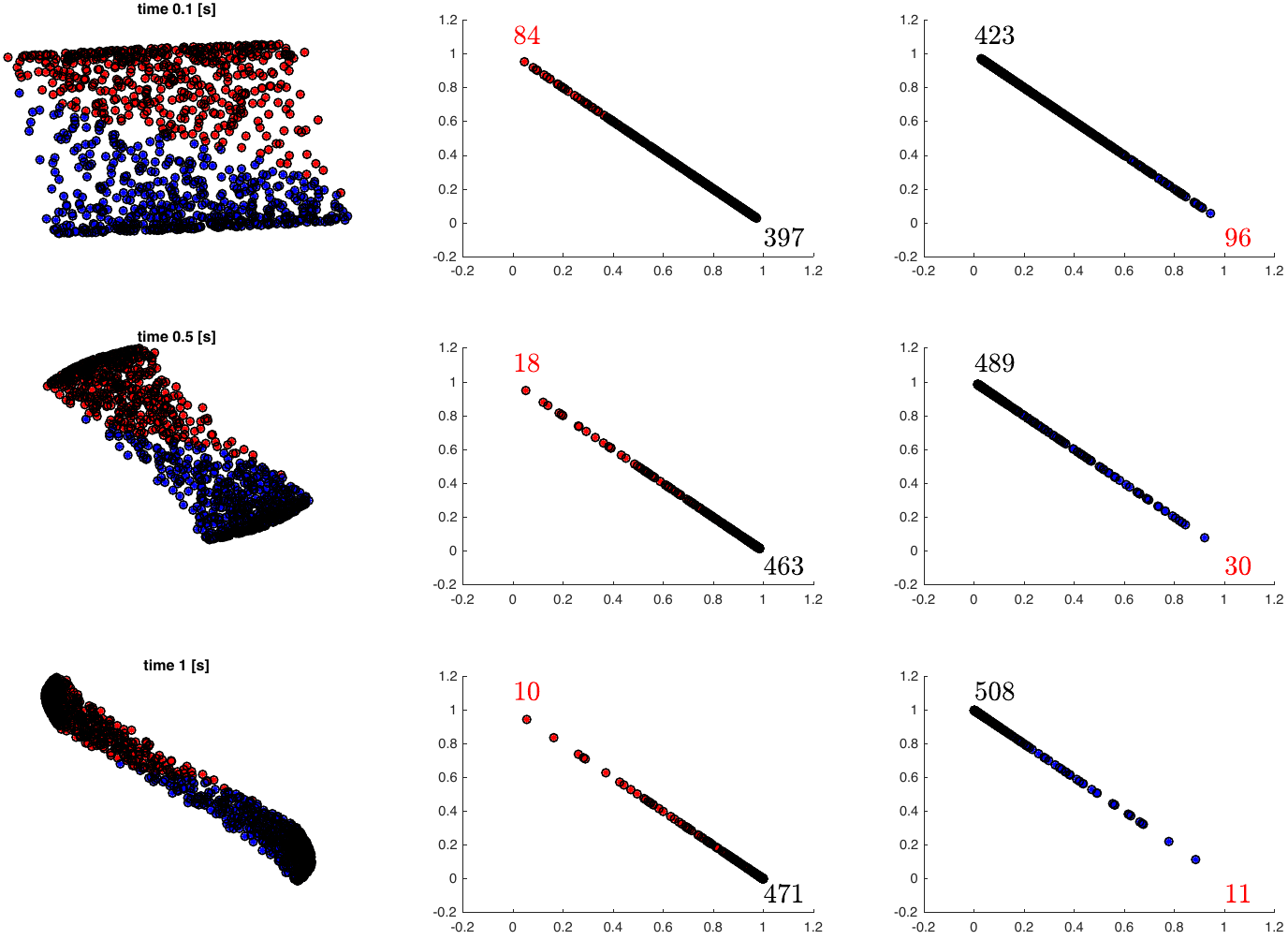}\\
  (a) & (b) 
	\end{tabular}
 \end{center}
  \caption{Graphical illustration of the classification results for: (a) Experiment 2; and (b) Experiment 3. The plots illustrate the connections between the temporal evolution of the ODE RNN $\tilde{\xi}^\theta(a(t))$ and the corresponding replicator dynamics for classification problems in the simplex.}
  \label{fig:experiment_2_3}
\end{figure*}

In the first experiment, we study the classification of an input $x\in\mathbb{R}^2$ into three possible labels $\mathcal{L}:=\{1,2,3\}$. The left plot in Figure \ref{fig:experiment_data} displays the training data set with the true labels depicted in green, red, and blue. The training data set consists of $D=3000$ points. In the second and third experiments, we consider the classification of an input in $x\in\mathbb{R}^2$ into two possible labels $\mathcal{L}:=\{1,2\}$. In this case, the center and right plots of Figure \ref{fig:experiment_data} showcase the training dataset with the true labels represented in red and blue. For this scenario, we used a training dataset with $D=1000$ points. 

For Experiment 1, even though the input data is in $\mathbb{R}^2$, the RNN lifts the dimensionality to $d_0=d_L=3$, as shown in the first and fourth column in Figure \ref{fig:experiment_1}. In this Figure, we also show the evolution over time of the population states in the simplex, for each individual label. The red numbers in the vertices of the simplex indicate the misclassified data elements, i.e., after training, only $37$ data points remain misclassified out of a total of $3000$ data points. Similarly, the results for Experiments 2 and 3 are shown in Figure \ref{fig:experiment_2_3}. In this figure, the evolution of the population state under the replicator dynamics is shown in the second, third, fifth, and sixth columns. After training, for both experiments we obtain that $21$ data points remain misclassified out of $1000$ data points.  It can also be observed in Figures \ref{fig:experiment_1} and \ref{fig:experiment_2_3} that the trajectories of the system approximate the vertices of the simplex as time increases, indicating a correct classification of the data with high probability. However, we note that for non-classification problems studied via RNNs with output activation functions given by the Boltzmann distribution, pursuing convergence towards an equilibrium point in the interior of the simplex might be of interest in the future via, e.g., evenhanded distributions of resources and/or errors as in \cite{BaTe_2018}. Finally, we note that the equivalence between the trajectories of the replicator system and the RNN ODE is preserved independently of the dimensionality of the input $x$ and the set of labels. However, a low input/output dimension was considered in this example to graphically illustrate the behavior of the RNN ODE.
\subsection{Classification of Digits using 10 Labels}

\begin{figure*}[t!]
   \centering
   \resizebox{0.95\textwidth}{!}{
   \begin{tabular}{c}
   \textcolor{black}{Input data samples $x \in [0,1]^{28 \times 28}$ to test the trained RNN}\\
   \includegraphics[clip,trim = 0cm 10cm 0cm 0.5cm,scale=0.8]{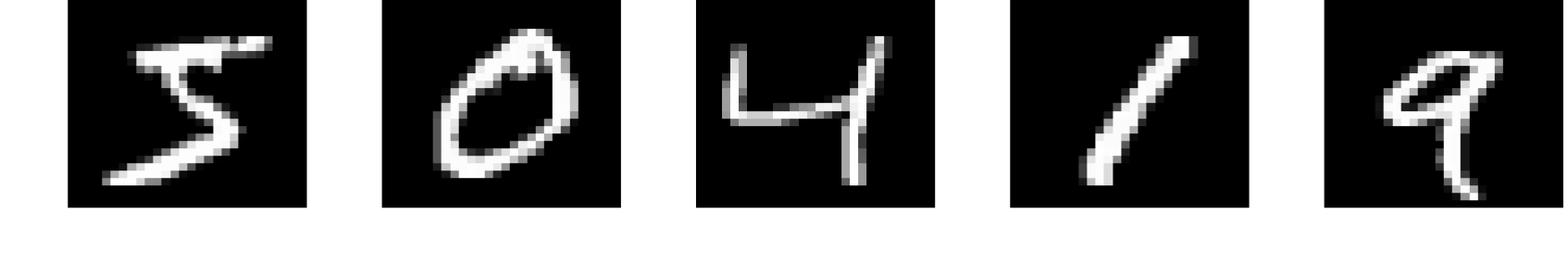}  \\
   \textcolor{black}{Payoff-related state $a(t)$, and population states (distribution) $y(t)$ at time $t=0.1$ [s] }\\
        \includegraphics[clip,trim = 0cm 0cm 0cm 5.5cm,scale=0.8]{Figs_2/DF_sample_1.pdf}  \\
        \\
        \textcolor{black}{Payoff-related state $a(t)$, and population states (distribution) $y(t)$ at time $5$ [s] }\\
        \includegraphics[clip,trim = 0cm 0cm 0cm 5.5cm,scale=0.8]{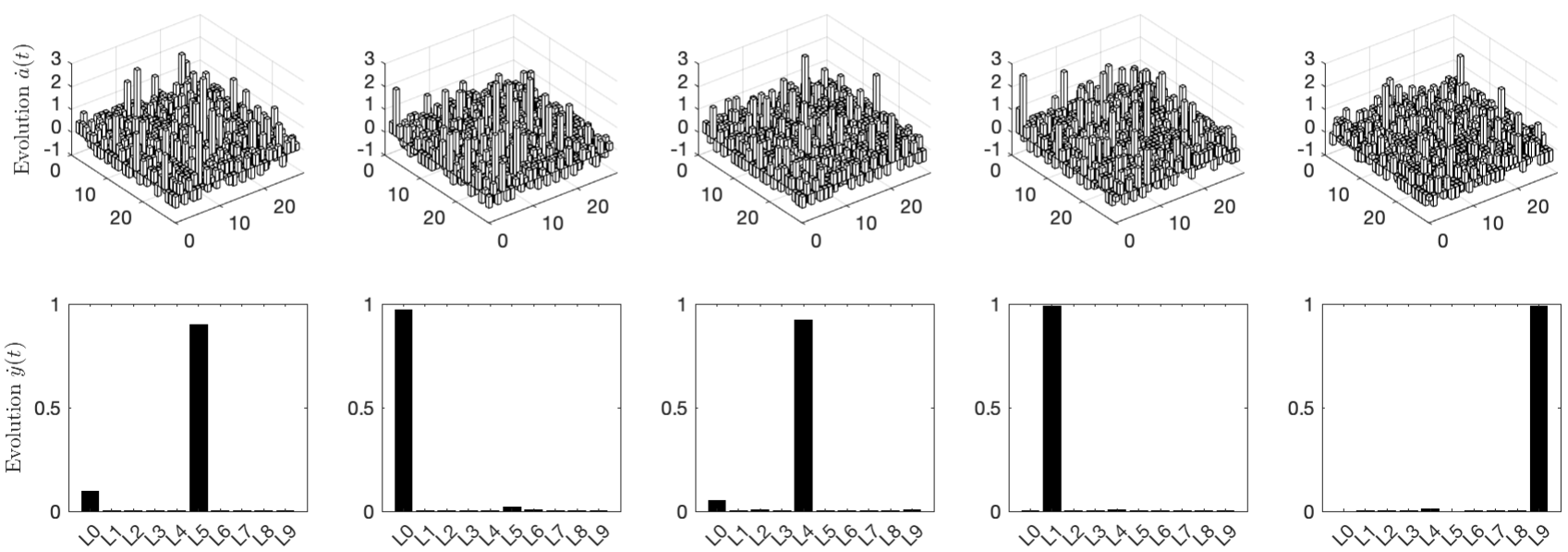}
   \end{tabular}
   }
   \caption{Classification of digits from black and white images. Initially, at the time $t=0.1$s, the classification of the sample figures is given by $(2,2,2,1,9)$ corresponding to the labels $(\mathrm{wrong},\mathrm{wrong},\mathrm{wrong},\mathrm{correct},\mathrm{correct})$. However, as shown in the last row of the figure, at the time $t = 5$s the classification is correct with high probability, i.e., the population states are located near the (correct) vertices of the simplex.}
   \label{fig:example}
\end{figure*}

\begin{figure}[h!]
   \centering
   \includegraphics[scale=0.37]{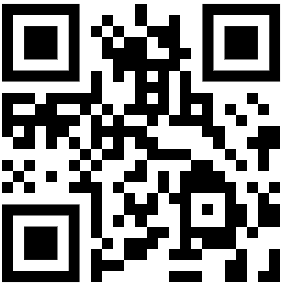}
      \vspace{-0.2cm}
   \caption{QR code for animated simulation of the behavior shown in Figure \ref{fig:example}.}
   \label{fig:video}
   \vspace{-0.4cm}
\end{figure}

Next, and to illustrate our results in a more practical problem, we consider the problem of classifying digits (MNIST) from black and white images. In this case, the dimension of the input data is substantially higher, namely, $x \in [0,1]^{784}$ for black and white images of dimension $28 \times 28$ pixels. We consider ten labels, one for each digit, such that $\mathcal{L}=\{1,2,\ldots,10\}$. In this case, while the simplex cannot be graphically illustrated as in Figures \ref{fig:scheme1} and \ref{fig:experiment_1}, we can still observe the evolution over time of the payoff dynamics $\dot{a}$ and the population dynamics $\dot{y}$ as a histogram, presented in Figure \ref{fig:example}. The selected architecture for the RNN satisfies 
$\phi_{\mathrm{in}}: [0,1]^{784} \to \mathbb{R}^{784}$, $\xi_\theta: \mathbb{R}^{784} \to \mathbb{R}^{784}$,  $\phi_{\mathrm{out}}: \mathbb{R}^{7840} \to [0,1]^{10}$, such that the evolution of $a \in \mathbb{R}^{28 \times 28}$ is obtained via \eqref{adynamics1} and the evolution of $y \in [0,1]^{10}$ is obtained via equation \eqref{firstRDs}. As observed in the last row of Figure \ref{fig:example}, in this case the classification problem is also correctly solved after 5 seconds. An animated simulation of the behavior of the system can also be observed via Figure \ref{fig:video}. 


\section{Conclusions and Future Research}
\label{sec:conclusions_nn}
A novel connection between RNNs and game-theoretic evolutionary dynamics was established by demonstrating that the output of an RNN ODE for classification problems follows a replicator system. This observation opens the door to studying the emerging behavior of RNNs using game-theoretic approaches and tools developed for population games, including passivity and dissipativity techniques, graphical models, and Lyapunov-based tools. The emerging feedback interconnection also suggests potential avenues for influencing the classification properties of the system via dynamic mechanism design. Future work will focus on these research directions.

\bibliographystyle{ieeetr}
\bibliography{references}

\begin{thebibliography}{10}

\bibitem{jafarpour2023efficient}
S.~Jafarpour, A.~Harapanahalli, and S.~Coogan, ``Efficient interaction-aware
  interval analysis of neural network feedback loops,'' {\em arXiv preprint
  arXiv:2307.14938}, 2023.

\bibitem{Chen2018}
R.~T.~Q. Chen, Y.~Rubanova, J.~Bettencourt, and D.~Duvenaud, ``Neural ordinary
  differential equations,'' in {\em 32nd Conference on Neural Information
  Processing Systems (NeurIPS)}, pp.~6571--6583, 2018.

\bibitem{krstic2023neural}
M.~Krstic, L.~Bhan, and Y.~Shi, ``Neural operators of backstepping controller
  and observer gain functions for reaction-diffusion {PDE}s,'' {\em arXiv
  preprint arXiv:2303.10506}, 2023.

\bibitem{liu2021second}
G.-H. Liu, T.~Chen, and E.~Theodorou, ``Second-order neural {ODE} optimizer,''
  {\em Ad. in Neural Inf. Proc. Sys.}, pp.~25267--25279, 2021.

\bibitem{marden2018game}
J.~R. Marden and J.~S. Shamma, ``Game theory and control,'' {\em Annual Rev. of
  Control, Robotics, and Aut. Syst.}, vol.~1, pp.~105--134, 2018.

\bibitem{brown2017studies}
P.~N. Brown and J.~R. Marden, ``Studies on robust social influence mechanisms:
  Incentives for efficient network routing in uncertain settings,'' {\em IEEE
  Ctrl. Syst. Magazine}, vol.~37, no.~1, pp.~98--115, 2017.

\bibitem{bhatia2014recurrent}
S.~Bhatia and R.~Golman, ``A recurrent neural network for game theoretic
  decision making,'' in {\em Proceedings of the Annual Meeting of the Cognitive
  Science Society}, vol.~36, 2014.

\bibitem{Hu2023}
R.~Hu and M.~Lauriere, {\em Recent Developments in Machine Learning Methods for
  Stochastic Control and Games}.
\newblock hal-03656245, 2023.

\bibitem{Effati2013}
S.~Effati and M.~Pakdaman, ``Optimal control problem via neural networks,''
  {\em Neural Computing and Applications}, vol.~23, no.~2013, pp.~2093--2100,
  2013.

\bibitem{Gomes2023b}
D.~Gomes, J.~Gutierrez, and M.~Lauriere, ``Machine learning architectures for
  price formation models,'' {\em Applied Mathematics and Optimization},
  vol.~88, no.~23, pp.~1--41, 2023.

\bibitem{BaChoBou2022}
J.~Barreiro-Gomez, S.~E. Choutri, and B.~Djehiche, ``Stability via adversarial
  training of neural network stochastic control of mean-field type,'' in {\em
  Proceedings of the 61st IEEE Control Conference on Decision and Control
  (CDC)}, pp.~7547--7552, 2022.

\bibitem{Vesseron2021}
N.~Vesseron, I.~Redko, and C.~Laclau, ``Deep neural networks are congestion
  games: From loss landscape to {Wardrop} equilibrium and beyond,'' in {\em
  International Conference in Artificial Intelligence and Statistics},
  pp.~1765--1773, 2021.

\bibitem{Ren2021}
C.~Ren, Z.~wu, D.~Xu, and W.~Xu, ``A game-theoretic analysis of deep naural
  networks,'' in {\em International Conference on Algorithmic Applications in
  Management}, pp.~369--379, 2021.

\bibitem{Tembine2019}
H.~Tembine, ``Deep learning meets game theory: Bregman-based algorithms for
  interactive deep generative adversarial networks,'' {\em IEEE Transactions on
  Cybernetics}, vol.~50, no.~3, pp.~1132--1145, 2020.

\bibitem{Jin2020}
Y.~Jin, Y.~Wang, M.~Long, J.~Wang, P.~S. Yu, and J.~Sun, ``A multi-player
  minimax game for generative adversarial networks,'' in {\em IEEE Int. Conf.
  on Multimedia and Expo}, pp.~1--6, 2020.

\bibitem{Chivukula2017}
A.~S. Chivukula and W.~Liu, ``Adversarial learning games with deep learning
  models,'' in {\em 2017 International Joint Conference on Neural Networks
  (IJCNN)}, pp.~2758--2767, 2017.

\bibitem{sandholm2010population}
W.~H. Sandholm, {\em Population games and evolutionary dynamics}.
\newblock MIT press, 2010.

\bibitem{hofbauer1998evolutionary}
J.~Hofbauer and K.~Sigmund, {\em Evolutionary games and population dynamics}.
\newblock Cambridge university press, 1998.

\bibitem{park2019population}
S.~Park, N.~C. Martins, and J.~S. Shamma, ``From population games to payoff
  dynamics models: A passivity-based approach,'' in {\em Conf. on Decision and
  Control}, pp.~6584--6601, 2019.

\bibitem{quijano2017role}
N.~Quijano, C.~Ocampo-Martinez, J.~Barreiro-Gomez, G.~Obando, A.~Pantoja, and
  E.~Mojica-Nava, ``The role of population games and evolutionary dynamics in
  distributed control systems: The advantages of evolutionary game theory,''
  {\em IEEE Control Systems Magazine}, vol.~37, no.~1, pp.~70--97, 2017.

\bibitem{BaTe_2018}
J.~Barreiro-Gomez and H.~Tembine, ``Constrained evolutionary games by using a
  mixture of imitation dynamics,'' {\em Automatica}, vol.~97, no.~2018,
  pp.~254--262, 2018.

\bibitem{arcak2020dissipativity}
M.~Arcak and N.~C. Martins, ``Dissipativity tools for convergence to nash
  equilibria in population games,'' {\em IEEE Transactions on Control of
  Network Systems}, vol.~8, no.~1, pp.~39--50, 2020.

\bibitem{pantoja2011population}
A.~Pantoja and N.~Quijano, ``A population dynamics approach for the dispatch of
  distributed generators,'' {\em IEEE Transactions on Industrial Electronics},
  vol.~58, no.~10, pp.~4559--4567, 2011.

\bibitem{BaObQu2017}
J.~Barreiro-Gomez, G.~Obando, and N.~Quijano, ``Distributed population
  dynamics: Optimization and control applications,'' {\em IEEE Trans. on Syst.,
  Man, and Cyb: Syst}, vol.~47, no.~2, pp.~304--314, 2017.

\bibitem{poveda2015shahshahani}
J.~I. Poveda and N.~Quijano, ``Shahshahani gradient-like extremum seeking,''
  {\em Automatica}, vol.~58, pp.~51--59, 2015.

\bibitem{sanfelice2010dynamical}
R.~S. G and A.~Teel, ``Dynamical properties of hybrid systems simulators,''
  {\em Automatica}, vol.~46, no.~2, pp.~239--248, 2010.

\end{thebibliography}
\end{document}